\documentclass[11pt,intlimits, draft,oneside]{amsart}

\usepackage[latin2]{inputenc}
\usepackage{latexsym}
\usepackage{amssymb}
\usepackage{amsthm}
\usepackage{amsmath}
\usepackage{enumerate}

\newtheorem{thm}{Theorem}[section]
\newtheorem{prop}[thm]{Proposition}
\newtheorem{lemma}[thm]{Lemma}

\theoremstyle{definition}
\newtheorem{defi}[thm]{Definition}
\newtheorem{remark}[thm]{Remark}
\newtheorem{history-remark}[thm]{Historical remark}
\newtheorem{example}[thm]{Example}
\newtheorem{problem}[thm]{Problem}
\newtheorem{question}[thm]{Question}

\newcommand\Eins{\mathbf{1}}

\renewcommand{\Re}{\mathrm{Re}\,}

\def\Fix{\mathop{\mathrm{Fix}}}
\def\rg{\mathop{\mathrm{rg}}}
\def\dd{\,d}

\def\clin{\mathop{\overline{\mathrm{lin}}}}

\def\Psigma{P\sigma}
\def\LLL{\mathcal{L}}

\def\ee{\mathrm{e}}
\def\RR{\mathbb{R}}
\def\NN{\mathbb{N}}
\def\CC{\mathbb{C}}
\def\lm{\rightarrow}

\def\la{\langle}
\def\ra{\rangle}

\begin{document}

\title{Weakly and almost weakly stable $C_0$-semigroups}

\author{Tanja Eisner}

\author{B\'{a}lint Farkas}

\address{B\'{a}lint Farkas \newline
AG Angewandte Analysis, Fachbereich Mathematik, Technische
Universit\"{a}t Darmstadt\newline Schlo{\ss}gartenstra{\ss}e 7,
D-64289 Darmstadt, Germany}
\email{farkas@mathematik.tu-darmstadt.de}

\author{Rainer Nagel}

\address{Tanja Eisner, Rainer Nagel \newline
Mathematisches Institut, Universit\"{a}t T\"{u}bingen\newline Auf
der Morgenstelle 10, D-70176, T\"{u}bingen, Germany}
\email{talo@fa.uni-tuebingen.de, rana@fa.uni-tuebingen.de}

\author{Andr\'{a}s Ser\'{e}ny}

\let\oldthefootnote\thefootnote\def\thefootnote{}
\footnotetext{The fourth author has been supported by the Marie
Curie Host Fellowship ``Spectral theory for evolution equations'',
contract number HPMT-CT-2001-00315.}
\let\thefootnote\oldthefootnote
\address{Andr\'{a}s Ser\'{e}ny \newline
Department of Mathematics and its Applications, Central European
University,\newline N\'{a}dor utca  9, H-1051 Budapest, Hungary }
\email{sandris@elte.hu}


\begin{abstract}
In this article we survey results concerning asymptotic properties
of $C_0$-semigroups on Banach spaces with respect to the weak
operator topology. The property ``no eigenvalues of
the generator on the imaginary axis'' is equivalent to weak
stability for 
most time values; a phenomenon called
``almost weak stability''. Further, sufficient conditions implying
the actual weak stability of a $C_0$-semigroup are also given. The
individual stability versions are considered as well, using both
boundedness of the local resolvent and integrability conditions for
the resolvent operator. By a number of examples we explain
weak and almost weak stability and illustrate the fundamental
difference between them. Where adequate, connections to strong
stability are observed, while many historical and
bibliographical remarks 
position the material in the
literature. We conclude the paper with some open questions and
comments.
\end{abstract}

\maketitle


\section{Introduction}

Strongly continuous semigroups on Banach spaces ($C_0$-semigroups
for short) provide a very efficient and elegant tool for the
treatment of concrete and abstract Cauchy problems (see Engel, Nagel
\cite{engel/nagel:2000}). They not only yield well-posedness results
(through the classical Hille-Yosida theorem and its variants), but
also allow a detailed description of important qualitative
properties of the solutions of the Cauchy problem. In this context,
it is 
important to describe the asymptotic behaviour of the
solutions. In terms of the semigroup, this means the following.
\begin{problem}\label{problem1}
Let $(T(t))_{t\geq 0}$ be a $C_0$-semigroup with generator $A$ on a
Banach space $X$. For each $x\in X$ describe the behaviour of
$T(t)x$ as $t\to\infty$.
\end{problem}
%
This question has been studied systematically in the monograph \cite{vanneerven:1996} by 
van Neerven, while Section V in Engel, Nagel \cite{engel/nagel:2000} 
or Sections V and VI in Engel, Nagel \cite{engel/nagel:2006} contain the major abstract 
results and some concrete applications.

For the following we assume $X$ to be reflexive and $(T(t))_{t\geq 0}$ to be bounded.
Then it follows from the Jacobs-Glicksberg-de Leeuw theory (see
\cite{engel/nagel:2000}, Thm.~V.2.8) that $X$ is the direct sum of
the two subspaces
\begin{align*}
X_r :=& \clin \bigl\{x \,\in\, D(A):\ Ax= i\lambda x \text{ for some } \lambda \in \RR \bigr\} \\
     =& \clin \bigl\{x\in X:\ T(t)x=\ee^{i\lambda t}x  \text{ for some } \lambda \in \RR  \bigr\},\\
X_s :=& \bigl\{ x \in X : 0 \text{ is a weak accumulation point of }
\{ T(t)x:\: t\geq 0 \} \bigr\}.
\end{align*}
The restriction of $(T(t))_{t\geq 0}$ to $X_r$ acts as a group and
its behaviour is determined by the purely imaginary eigenvalues
$i\lambda \in \Psigma(A)$ of the generator $A$ (as usual
$\Psigma(A)$ denotes the point spectrum of the operator $A$), or,
equivalently, by the periodic orbits of $T(t)$. Therefore, for this
part of $(T(t))_{t\geq 0}$ Problem \ref{problem1} is solved in a
reasonable way.

On the other hand, $A$ has no eigenvectors with purely imaginary
eigenvalue in $X_s$  and, based on the above characterisation, one
expects 
\begin{equation*}
``\lim_{t\to\infty}T(t)x=0"   \ \ \text{ for all } x\in X_s
\end{equation*}
in some sense. Therefore, Problem \ref{problem1} reduces to the
following.
\begin{problem}\label{problem2}
Let $(T(t))_{t\geq 0}$ be a bounded $C_0$-semigroup with generator
$A$ on a reflexive Banach space $X$. Find necessary and sufficient
conditions for $(T(t))_{t\geq 0}$ to be \emph{stable}, i.e., to
satisfy
\begin{equation*}
\lim_{t\to\infty}T(t)x=0 \ \ \text{ for all } x\in X_s,
\end{equation*}
where the limit is taken in some appropriate topology.
\end{problem}

The first result in this direction goes back to Lyapunov.
In 1892 he gave the following characterisation for matrix semigroups.
\begin{thm}
Let $(T(t))_{t\geq 0}$ be a $C_0$-semigroup with generator $A$ on a
finite dimensional Banach space $X$. The semigroup $(T(t))_{t\geq
0}$ is uniformly stable, i.e., $\lim\limits_{t\to\infty}\| T(t)\|
=0$ if and only if
\begin{equation*}
s(A):=\sup\big\{\Re \lambda:\ \lambda \in\sigma(A)\bigr\}<0.
\end{equation*}
\end{thm}
This characterisation does no longer hold on infinite dimensional
spaces and for  unbounded generators (see \cite{engel/nagel:2000},
Sect.~IV.3 for counterexamples). Still, the situation concerning uniform stability is quite well
understood and there are many results, using additional properties
of the Banach space and/or the semigroup, on uniformly stable
semigroups.  We refer to Engel, Nagel \cite{engel/nagel:2000}, Sect.~V.1.b, Arendt, Batty, Hieber, Neubrander \cite{arendt/etal:2001}, Sect.~5.2-5.3, and van Neerven \cite{vanneerven:1996}, Ch.~3. 

The \emph{strong stability} of $(T(t))_{t\geq 0}$, i.e., the
property that
\begin{equation*}
\lim_{t\to\infty}\| T(t)x \| =0 \ \ \text{ for all } x\in X,
\end{equation*}
is not so well understood and only in 1988 did Arendt, Batty
\cite{arendt/batty:1988}, Lyubich, V\~u \cite{lyubich/vu:1988}
obtain the following sufficient, but not necessary condition.
\begin{thm}
Let $(T(t))_{t\geq 0}$ be a bounded $C_0$-semigroup with generator
$A$ on a reflexive Banach space $X$ with $P_\sigma(A)\cap i\RR =\emptyset$. If $\sigma(A)\cap i\RR$ is
countable, then $(T(t))_{t\geq 0}$ is strongly stable.
\end{thm}
More results on strong stability with even a necessary and
sufficient condition in the case of Hilbert spaces are due to
Tomilov \cite{tomilov:2001} and Chill, Tomilov
\cite{chill/tomilov:2003, chill/tomilov:2004}.

Surprisingly, compared to the above there are relatively few results
on \emph{weakly stable} semigroups $(T(t))_{t\geq 0}$, i.e.,
semigroups satisfying
\begin{equation*}
\lim_{t\to\infty} \langle T(t)x,\varphi \rangle = 0 \ \ \text{ for
all } x\in X \text{ and all } \varphi \in X'.
\end{equation*}

Note that, for example, isometric semigroups (in particular unitary groups)
never converge strongly (except in the trivial case), so convergence
with respect to the weak topology is the natural mode of convergence
for such semigroups.


The absence of a theory for weak stability of semigroups   
is regrettable not just for pure mathematical reasons. In the
setting of quantum theory, one thinks of $X$ as the state space,
while $X'$ is the space of observables of some system. Therefore,
the scalar valued function
\begin{equation*}
t\mapsto  \langle T(t)x,\varphi \rangle
\end{equation*}
gives the time evolution of a measuring process. Consequently, weak
stability is the property of a system which can indeed be observed.

Also in ergodic theory, the concept of weak stability occurs
naturally (under the name ``strong mixing''). We quote here the
following from Katok, Hasselblatt \cite{katok/hasselblatt:1995},
p.~748.

\begin{quote}\it
``... It [strong mixing] is, however, one of those notions, that is
easy and natural to define but very difficult to study...''
\end{quote}
We will see later that a weaker concept (``almost weak stability'', see 
Section \ref{section:characterization})
is in fact relatively
easy to characterise, whereas our knowledge about weak stability
itself is still limited.


\smallskip

 The aim of this paper is thus to survey the known
results on weak asymptotic properties of $C_0$-semigroups, and also
to propose some ideas for further research.

We will start with the essentially easier notion of
 \emph{``almost weak stability''}. By virtue of the
Jacobs-Glicksberg-de Leeuw decomposition, for bo\-unded semigroups
on reflexive spaces, almost weak stability turns out to be
equivalent to
\begin{center}\it
the generator $A$ has no purely imaginary eigenvalues.
\end{center}
 In Section
\ref{section:characterization} we will give a series of equivalent
properties justifying the term ``almost weakly stable'', and give a
brief account on the history of this and related notions. In
particular, we mention \emph{weakly almost periodic functions} to
which we return in the concluding Section
\ref{section:comments_and_open_quetions} as well.

\smallskip

In Section \ref{section:weak_stability}, we proceed to weak
stability. In Hilbert spaces, one can decompose contraction
semigroups into a weakly stable and a unitary part, which allows us
to restrict the study to unitary groups. In general Banach spaces,
we present sufficient conditions for convergence to $0$ of orbits of
$C_0$-semigroups in terms of integrability of the resolvent of the
generator. We also discuss how the convergence of the semigroups
$(T(t))_{t\geq 0}$ and $(T(n))_{n\in\NN}$ is related. 

\smallskip

We devote Section \ref{section:examples} to examples and
counterexamples to illustrate the notions of weak and almost weak
stability. In particular, we show that in a certain sense only a
minority of almost weakly stable semigroups are indeed weakly stable.

\smallskip

Section \ref{section:individual_stability} is devoted to the property 
\begin{equation*}
\lim_{t\to\infty} \langle T(t)x, y \rangle =0 \quad \text{ for some given  } x \text {  and } y, 
\end{equation*}
the so-called \emph{weak individual stability}. 
%
We investigate this phenomenon assuming the existence of
a bounded local resolvent. Here the borderline between strong and
weak individual stability is very narrow: in the presence of certain
geometric properties of the underlying Banach space one even obtains
strong stability. We also elaborate on this aspect by recalling
several results from the literature.

\smallskip

The last Section \ref{section:comments_and_open_quetions} is
intended to reveal connections to results which do not exactly fit
into the line of this survey, but are connected to our topic and
thus deserve a few words.  We finally pose some open questions.

\section{Almost  weak stability} \label{section:characterization}

Let us start by discussing almost weak stability of
$C_0$-semigroups, a concept which is close to weak stability but
much easier to investigate.

Our main functional analytic tool will be (relative) compactness for
the weak operator topology. We denote by $\LLL_\sigma(X)$ the space
$\LLL(X)$ endowed with the weak operator topology and recall the
following characterisation of relative compactness in
$\LLL_\sigma(X)$.

\begin{lemma}[{\cite{engel/nagel:2000}, Lemma V.2.7}] \label{weakcomp}
For a set of operators $\mathcal{T}\subset \LLL(X)$, $X$ a Banach
space, the following assertions are equivalent.
\begin{enumerate}[(a)]
\item $\mathcal{T}$ is relatively compact in $\LLL_\sigma(X)$.
\item $\{Tx:\:T\in \mathcal{T}\}$ is relatively weakly compact in $X$ for all $x\in X$.
\item $\mathcal{T}$ is bounded, and  $\{Tx:\:T\in \mathcal{T}\}$ is relatively weakly compact in $X$ for all $x$ in some dense subset of $X$.
\end{enumerate}
\end{lemma}

\noindent Let us give some examples of relatively weakly compact
subsets of operators.
\begin{example}\label{ex:relcomp} \rule{0pt}{0pt}
\begin{enumerate}[(a)]
\item On a reflexive Banach space $X$ any norm-bounded family $\mathcal{T}\subseteq \LLL(X)$ is relatively weakly compact.
\item Let $\mathcal{T}\subseteq \LLL(L^1(\mu))$ be a norm-bounded subset of positive operators on the Banach lattice $L^1(\mu)$, and suppose that $Tu\leq u$ for some $\mu$-almost everywhere positive $u\in L^1(\mu)$ and every $T\in \mathcal{T}$. Then $\mathcal{T}$ is relatively
weakly compact since the order interval $[-u,u]$ is weakly compact
and $\mathcal{T}$-invariant (see Schaefer \cite{schaefer:1974},
Thm.~ II.5.10 (f) and Prop.~II.8.3).
\item \label{ex:semitop} Let $S$ be a semi-topological semigroup, i.e., a (multiplicative) semigroup $S$ which is a
topological space such that the multiplication is separately
continuous (see Engel, Nagel \cite{engel/nagel:2000}, Sec.~V.2).
Consider the space $C(S)$ of bounded, continuous (real- or
complex-valued) functions over $S$.  For $s\in S$ define the
corresponding rotation operator $(L_s f)(t):= f(s\cdot t)$. A
function $f\in C(S)$ is said to be \emph{weakly almost periodic} if
the set $\{L_s f:\: s\in S\}$ is relatively weakly compact in
$C(S)$, see Berglund, Junghenn, Milnes \cite{berglund/etal:1989},
Def.~4.2.1 (cf also Section \ref{subsec:weakper}). The set of weakly
almost periodic functions is denoted by $WAP(S)$. If $S$ is a
\emph{compact} semi-topological semigroup, then $C(S)=WAP(S)$ holds,
see \cite{berglund/etal:1989}, Cor.~4.2.9. This means that for a
compact semi-topological semigroup $S$ the set $\{L_s:\:s\in S\}$ is
always relatively weakly compact. We come back to this example in
Example \ref{fluss} and in the proof of Theorem \ref{str.stable_1}.
\end{enumerate}
\end{example}
\noindent We now turn our attention to $C_0$-semigroups (see Engel,
Nagel \cite{engel/nagel:2000} for the general theory).

\begin{defi} A $C_0$-semigroup $(T(t))_{t\geq 0}$ on a Banach space $X$ is called \emph{relatively weakly
compact} if the set $\mathcal{T}:=\{T(t):\: t\geq 0\}$ satisfies one
of the equivalent conditions in Lemma \ref{weakcomp}.
\end{defi}
In the following we will concentrate only on relatively weakly
compact semigroups. We note that every weakly stable semigroup
has this property, therefore this is not a strong restriction
with respect to our aims.

The following property of relatively weakly compact semigroups will
be used in the sequel without explicit reference, see Engel, Nagel
\cite{engel/nagel:2000}, Sec.~V.4.
\begin{prop} Let $(T(t))_{t\geq 0}$ be a relatively weakly compact semigroup on a Banach space $X$. Then it is \emph{mean ergodic}, i.e., we have (weak and even strong) convergence of the Ces\`{a}ro means
\begin{equation*}
\lim_{t\to\infty}\frac{1}{t}\int_0^t T(s)x\dd s = Px\qquad\mbox{for
all $x\in X$},
\end{equation*}
where $P\in \LLL(X)$ is a projection onto $\Fix (T):=\bigcap_{t\geq
0}\Fix (T(t))$, the so-called \emph{ergodic projection}.
\end{prop}

Assume now $(T(t))_{t\geq 0}$ to be relatively weakly compact. By
the decomposition theorem of Jacobs-Glicksberg-de Leeuw the Banach
space $X$ is a direct sum of the two invariant subspaces
\begin{align*}
X_r =& \clin \bigl\{x \,\in\, D(A):\ Ax= i\lambda x \text{ for some } \lambda \in \RR \bigr\}, \\
X_s =& \bigl\{ x \in X : 0 \text{ is a weak accumulation point of }
\{ T(t)x:\: t\geq 0 \} \bigr\}.
\end{align*}
(See Maak \cite{maak:1954}, Jacobs \cite{jacobs:1955}, Glicksberg,
de Leeuw \cite{glicksberg/deleeuw:1961} and ultimately Krengel
\cite{krengel:1985}, Sec.~2.4 for detailed discussion and historical
remarks and also Engel, Nagel \cite{engel/nagel:2000}, Thm.~V.2.8).

So we see that the property ``no eigenvalues of the generator on the imaginary axis'' is equivalent to the fact that  $0$ is a weak accumulation point of the orbit $\{T(t)x:\ t\geq 0\}$ for every $x\in X$. The following theorem gives more detailed information about the asymptotic behaviour of the orbits in this case.

%
%

\begin{thm}\label{str.stable_1}\rule{0pt}{0pt}
Let  $(T(t))_{t \geq 0}$ be a relatively weakly compact
$C_0$-semigroup on a Banach space $X$ with generator $A$.
The following assertions are equivalent.
\begin{enumerate}[(i)]
%
\item[(i)] $\ 0\in \overline{\{T(t)x:\:t\geq 0\}}^{\sigma}$ for every $x\in X$;

\item[(i')] $\ 0\in \overline{\{T(t):\:t\geq 0\}}^{\LLL_\sigma}$;

\item[(ii)]  For every $x\in X$ there exists a sequence $\{t_n\}_{n=1}^\infty$
with $t_n\to \infty$ such that $T(t_n)x \stackrel{\sigma}{\to} 0$;

\item[(iii)]  For every $x\in X$ there exists a set $M\subset \RR_+$ with
density $1$ such that $T(t)x\stackrel{\sigma}{\to} 0$, as $t\in
M,\:t\to \infty$;
\item[(iv)] $\displaystyle \tfrac{1}{t}\smallint_0^t |\langle T(s)x,y\rangle | \dd s \underset{{t\to \infty}}{\longrightarrow} 0$
for all $x\in X$, $y\in X'$;
\item[(v)] $\displaystyle \lim_{ a \to 0+} a \smallint_{-\infty}^\infty | \langle R(a+i s,A)x,y \rangle |^2
\dd s = 0$  for all $x\in X$, $y\in X'$;
\item[(vi)] $\displaystyle \lim_{a\to 0+} a R(a+i s,A)x = 0$ for all $x\in X$ and $s\in \RR$;
\item[(vii)] $\Psigma (A)\cap i\RR = \emptyset$, i.e., $A$ has no purely imaginary eigenvalues.
\end{enumerate}
If, in addition, $X'$ is separable, then the conditions above are
also equivalent to
\begin{enumerate}[(i)]
\item[(ii$^*$)] There exists a sequence $\{t_n\}_{n=1}^\infty$ with $t_n\to \infty$ such that
$T(t_n)\stackrel{\sigma}{\to} 0$;
\item[(iii$^*$)] There exists a set $M\subset \RR_+$ with density $1$ such that
$T(t)\stackrel{\sigma}{\to}0$, $t\in M$ and $t\to\infty$.
\end{enumerate}
\end{thm}
Recall that the (asymptotic) density of a measurable set $M\subset
\RR_+$ is
\begin{equation*}d(M):=\lim_{t\to \infty} \frac{1}{t} \lambda ([0,t]\cap M),\end{equation*}
whenever the limit exists (here $\lambda$ is the Lebesgue measure on
$\RR$). 
%

We will use the following classical lemma (for the proof see
Petersen \cite{petersen:1983}, Lemma 6.2).
\begin{lemma}\label{dens1} Let $f:\RR_+ \to \RR_+$ be continuous and bounded. The following
assertions are equivalent.
\begin{enumerate}[(a)]
\item $\displaystyle\tfrac{1}{t}\smallint_0^t f(s)\dd s \to 0$ as $ t\to \infty$;
\item There exists a set $M\subset \RR_+$ with density $1$ such that $f(t)\to 0,
\ \ t\in M$ and $t\to \infty$.
\end{enumerate}
\end{lemma}

\noindent
\begin{proof}[Proof of Theorem \ref{str.stable_1}] \rule{0pt}{0pt}

The proof of the implication (i') $\Rightarrow$ (i) is trivial. 
The implication (i) $\Rightarrow$ (ii) holds since in Banach spaces
weak compactness and weak sequential compactness coincide (see
Eberlein-\v Smulian theorem, e.g., Dunford, Schwartz
\cite{dunford/schwartz:1958}, Thm.~V.6.1).

\smallskip\noindent If (vii) does not hold, then (ii) can not be true
by  the spectral mapping theorem for the point spectrum (see Engel,
Nagel \cite{engel/nagel:2000}, Thm.~V.3.7), hence (ii) $\Rightarrow$
(vii).

\smallskip\noindent The implication (vii) $\Rightarrow$ (i') is the
main consequence of the Jacobs-Glicksberg-de Leeuw theorem and
follows from the construction in its proof, see Engel, Nagel
\cite{engel/nagel:2000}, p.~313.

This proves the equivalences (i) $\Leftrightarrow$ (i')
$\Leftrightarrow$ (ii) $\Leftrightarrow$ (vii).

\smallskip\noindent (vi) $\Leftrightarrow$ (vii): Since the semigroup
$(T(t))_{t\geq 0}$ is mean ergodic and bounded, the decomposition
$X=\ker A \oplus \overline{\rg A}$ holds (see \cite{engel/nagel:2000}, Lemma
V.4.4). This implies by the mean ergodic theorem (see, e.g., Arendt,
Batty, Hieber, Neubrander \cite{arendt/etal:2001}, Cor.~4.3.2) that
the limit
\begin{equation*}
Px:=\lim_{a\to 0+} aR(a,A)x
\end{equation*}
exists for all $x\in X$ with a projection $P$ onto $\ker A$.
Therefore, $0\notin \Psigma (A)$ if and only if $P=0$. Take now
$s\in \RR$. The semigroup $(\ee^{i st}T(t))_{t\geq 0}$ is also
relatively weakly compact and hence mean ergodic. Repeating the
argument for this semigroup we obtain (vi) $\Leftrightarrow$ (vii).

\smallskip\noindent (i') $\Rightarrow$ (iii): Let
$S:=\overline{\{T(t):\:t\geq 0\} }^{\LLL_\sigma}\subseteq\LLL(X)$
which is a compact semi-topological semigroup if considered with the
usual multiplication and the weak operator topology. By (i) we have
$0\in S$. Define the operators $\tilde{T}(t):C(S) \to C(S)$ by
\begin{equation*}
(\tilde{T}(t)f)(R):= f(T(t)R), \ \ \ f\in C(S),\: R\in S.
\end{equation*}
By Nagel (ed.) \cite{nagel:1986}, Lemma B-II.3.2,
$(\tilde{T}(t))_{t\geq 0}$ is a $C_0$-semigroup on $C(S)$.

By Example \ref{ex:relcomp} (\ref{ex:semitop}) the set
$\{f(T(t)\:\cdot):\: t\geq 0\}$ is relatively weakly compact in
$C(S)$ for every $f\in C(S)$. It means that every orbit
$\{\tilde{T}(t)f:\: t\geq 0\}$ is relatively weakly compact, and,
 by Lemma \ref{weakcomp}, $(\tilde{T}(t))_{t\geq0}$ is a
relatively weakly compact semigroup.

\smallskip Denote by $\tilde{P}$ the mean ergodic projection of
$(\tilde{T}(t))_{t\geq0}$. We have $\Fix(\tilde{T})=\bigcap_{t\geq
0}\Fix(\tilde{T}(t))=\langle \Eins \rangle$. Indeed, for $f\in
\Fix(\tilde T)$ one has $f(T(t)I)=f(I)$ for all $t\geq 0$ and
therefore $f$ should be constant. Hence  $\tilde{P}f$ is constant
for every $f\in C(S)$. By definition of the ergodic projection
\begin{equation}\label{eq:meaner}
(\tilde{P}f)(0)=\lim_{t\to\infty} \frac{1}{t} \int_0^t
\tilde{T}(s)f(0)\dd s= f(0).
\end{equation}
Thus we have
\begin{equation}\label{eq:tildeP}
(\tilde{P}f)(R)= f(0)\cdot \Eins,\qquad f\in C(S),\: R\in S.
\end{equation}

\noindent Take now $x\in X$. By Theorem 3 and its proof in Dunford,
Schwartz \cite{dunford/schwartz:1958}, p.~434, the weak topology on
the orbit $\{T(t)x:\ t\geq 0\}$ is metrisable and coincides with the
topology induced by some sequence $\{y_n\}_{n=1}^\infty\subset
X'\setminus\{0\}$. Consider $f_{x,n}\in C(S)$ defined by
\begin{equation*}
f_{x,n}(R):= |\langle Rx, \tfrac{y_n}{\|y_n\|}\rangle|, \qquad R\in
S,
\end{equation*}
and $f_x\in C(S)$ defined by
\begin{equation*}
f_x(R):= \sum_{n\in\NN}\frac{1}{2^n}f_{x,n}(R), \qquad R\in S.
\end{equation*}
By \eqref{eq:tildeP} we obtain
\begin{equation*}
0= \lim_{t\to \infty} \frac{1}{t} \int_0^t \tilde{T}(s)f_{x,y}
(I)\dd s = \lim_{t\to \infty} \frac{1}{t} \int_0^t f_x(T(s)) \dd s.
\end{equation*}
Lemma \ref{dens1} applied to the continuous and bounded function
$\RR_+\ni t \mapsto f(T(t)I)$ yields a set $M\subset\RR$ with
density 1 such that
\begin{equation*}
f_x(T(t))\to 0 \quad \text{ as } t\to \infty,\ \ t\in M.
\end{equation*}
By definition of $f_x$ and by the fact that the weak topology on the
orbit is induced by  $\{y_n\}_{n=1}^\infty$  we have in particular
that
\begin{equation*}
T(t)x \stackrel{\sigma}{\to} 0 \quad \text{ as } t\to \infty,\ \
t\in M.
\end{equation*}
This proves (iii).

\smallskip\noindent (iii)  $\Rightarrow$ (iv) follows directly from
Lemma \ref{dens1}.

\smallskip\noindent (iv)  $\Rightarrow$ (vii) holds by  the spectral
mapping theorem for the point spectrum (see Engel, Nagel
\cite{engel/nagel:2000}, Thm.~V.3.7).

\smallskip\noindent (iv) $\Leftrightarrow$ (v): Clearly, the
semigroup $(T(t))_{t\geq 0}$ is bounded. Take $x\in X$, $y\in X'$
and let $a>0$. By the Plancherel theorem applied to the function
$t\mapsto \ee^{-at}\langle T(t)x,y\rangle$ we have
\begin{equation*}
\int_{-\infty}^\infty |\langle R(a+i s,A)x,y \rangle|^2 \dd s = 2\pi
\int_0^\infty \ee^{-2at}|\langle T(t)x,y\rangle|^2 \dd t.
\end{equation*}
We obtain by the equivalence of Abel and Ces\`{a}ro limits (see, e.g.,
Hardy \cite{hardy:1949}, p.~136)
\begin{eqnarray}
\lim_{a\to 0+} a \int_{-\infty}^\infty |\langle R(a+i s,A)x,y
\rangle|^2 \dd s
&=& 2\pi \lim_{a\to 0+} a \int_0^\infty \ee^{-2at}|\langle T(s)x,y\rangle|^2 \dd s \nonumber \\
&=& \pi \lim_{t\to \infty} \frac{1}{t} \int_0^t |\langle
T(s)x,y\rangle|^2 \dd s. \label{eq:plancherel}
\end{eqnarray}

\noindent Note that for a bounded continuous function $f:\RR_+ \to
\RR_+$ with $C:=\sup f(\RR_+)$ we have
\begin{equation*}
\biggl(\frac{1}{Ct} \int_0^t f^2(s)\dd s\biggr)^2 \leq \biggl(
\frac{1}{t} \int_0^t f(s)\dd s \biggr)^2 \leq \frac{1}{t} \int_0^t
f^2(s)\dd s,
\end{equation*}
which together with \eqref{eq:plancherel} gives the equivalence of
(iv) and (v).

\smallskip

For the additional part of the theorem suppose $X'$ to be separable.
Then so is $X$, and we can take dense subsets $\{x_n\neq 0:\
n\in\NN\}\subseteq X$ and $\{y_m\neq 0:\ m\in\NN\}\subseteq X'$.
Consider the functions
\begin{equation*}
 f_{n,m}:S\to\RR,\qquad f_{n,m}(R):=\bigl|\bigl\langle R \tfrac{x_n}{\|x_n\|},\tfrac{y_m}{\|y_m\|}\bigr\rangle \bigr|,\quad n,m\in \NN,
\end{equation*}
which are continuous and uniformly bounded in $n,m\in\NN$. Define
the function
\begin{equation*}
f:S\to\RR, \quad \ \
f(R):=\sum_{n,m\in\NN}\frac{1}{2^{n+m}}f_{n,m}(R).
\end{equation*}
Then clearly $f\in C(S)$. Thus, as in the proof of the implication
(i') $\Rightarrow$ (iv), i.e., using \eqref{eq:meaner} we obtain
\begin{equation*}
\frac{1}{t}\int_{0}^t f(T(s)I)\dd s
\underset{{t\to\infty}}{\longrightarrow} 0.
\end{equation*}
Hence, applying Lemma \ref{dens1} to the continuous and bounded
function $\RR_+\ni t \mapsto f(T(t)I)$ gives the existence of a set
$M$ with density $1$ such that $f(T(t))\to 0$ as $t\to \infty$,
$t\in M$. In particular, $|\langle T(t)x_n,y_m\rangle|\to 0$ for all
$n,m\in\NN$ as $t\in M$, $t\to\infty$, which, together with the
boundedness of $(T(t))_{t\geq 0}$, proves the implication (i')
$\Rightarrow$ (iii$^*$). The implications (iii$^*$) $\Rightarrow$
(ii$^*$) $\Rightarrow$ (ii') are straightforward, hence the proof is
complete.
\end{proof}

The above theorem shows that starting from ``no purely imaginary eigenvalues of the ge\-ne\-rator'', one arrives at properties like (iii) on the asymptotic behaviour of the orbits of the semigroup. This justifies the following name for this property. 
\begin{defi}
We will call a relatively weakly compact $C_0$-semigroup \emph{almost weakly stable} if it satisfies one of the equivalent conditions in Theorem \ref{str.stable_1}. 
\end{defi}
\begin{history-remark}
Theorem \ref{str.stable_1} and especially the implication (vii)
$\Rightarrow$ (iii) has a long history. It goes back to ergodic
theory and von Neumann's spectral mixing theorem for flows, see
Halmos \cite{halmos:1956}, Mixing Theorem, p.~39. This has been
generalised to operators on Banach spaces by many authors, see,
e.g.,  Nagel \cite{nagel:1974}, Jones, Lin \cite{jones/lin:1976, jones/lin:1980} and Krengel
\cite{krengel:1985}, pp.~108--110. The
implication (vii) $\Rightarrow$ (i) appears also in Ruess, Summers
\cite{ruess/summers:1992jmaa} (see also Section
\ref{subsec:weakper}). The conditions (i), (iii) and (iv) were
studied by Hiai \cite{hiai:1978} also for strongly measurable
semigroups. He related it to the discrete case as well. See also
K\"uhne \cite{kuehne:1981, kuehne:1987}.
\end{history-remark}

\begin{remark}
The conditions in Theorem \ref{str.stable_1} are of quite different
nature. Conditions (i)--(iv) as well as (ii$^*$) and (iii$^*$) give
information on the behaviour of the semigroup $(T(t))_{t \geq 0}$,
while conditions (v)--(vii) deal with the generator and its
resolvent. Among them condition (vii) apparently is the simplest to
verify.
\end{remark}
\begin{remark}\label{rem:strongfulint}
It is surprising that there is a characterisation of strong
stability of bounded $C_0$-semigroups on Hilbert spaces which is
completely analogous to the equivalence (i') $\Leftrightarrow$ (v)
in Theorem \ref{str.stable_1}. More precisely, a bounded
$C_0$-semigroup $(T(t))_{t \geq 0}$ on a Hilbert space $H$ with
generator $A$ is strongly stable if and only if
\begin{equation*} \lim_{a \to 0+} a \int_{-\infty}^\infty \|R(a+i s,A)x\|^2 \dd s = 0\end{equation*}
holds for every $x \in H$ (see Tomilov \cite{tomilov:2001},
pp.~108--110).
  We note that also the proofs of the equivalence (i') $\Leftrightarrow$ (v) and
Theorem~3.1 in \cite{tomilov:2001} are analogous.
\end{remark}

\begin{remark}
The conditions (iii) and (iii$^*$) show that all the orbits
$t\mapsto T(t)x$ converge weakly to $0$ as $t\rightarrow\infty$ for
$t$ in a large subset. 
In general, it may happen however
that this large set is not $\mathbb{R}_+$, i.e., $(T(t))_{t\geq 0}$
is not weakly stable (for examples see Section
\ref{section:examples}). Here is the essential difference to strong
stability: for a bounded semigroup $(T(t))_{t\geq 0}$ the
convergence $\|T(t_n)x\|\rightarrow 0$ for a sequence
$t_n\rightarrow \infty$ already implies $\displaystyle
\lim_{t\rightarrow\infty} \|T(t)x\| = 0$.
\end{remark}

\section{Weak stability} \label{section:weak_stability}

As we have already noted in the introduction, a bounded
$C_0$-semigroup on a reflexive Banach space $X$ induces the
Jacobs-Glicksberg-de Leeuw decomposition $X=X_r \oplus X_s$, but
orbits in $X_s$ do not necessarily converge (weakly) to zero.
However, in the particular case of contractive semigroups on Hilbert
spaces, one can detach the subspace of all weakly stable orbits and
characterise its complement.

\begin{thm}
Let $(T(t))_{t\geq0}$ be a $C_0$-semigroup of contractions on a
Hilbert space $H$ and define
\begin{equation*}W:=\bigl\{ x \in H: \ \lim_{t\to\infty} \langle T(t)x,x \rangle = 0 \, \bigr\} .\end{equation*}
Then $W$ is a closed subspace of $H$, $W$ and $W^{\perp}$ are
$(T(t))_{t\geq 0}$-invariant, the restricted semigroup
$(T(t)_{\mid_W})_{t\geq 0}$ is weakly stable on $W$ and
$(T(t)_{\mid_{W^{\perp}}})_{t\geq 0}$ is unitary on $W^{\perp}$.
\end{thm}

\noindent For the proof we refer the reader to Luo, Guo, Morgul
\cite{luo/gou/morgul:1999}, Thm.~3.18, p.~122, or see Foguel
\cite{foguel:1963} for the analogous discrete case.

In the following propositions we state some immediate consequences
of the above decomposition.

\begin{prop}
Let $(T(t))_{t\geq0}$ be a $C_0$-semigroup of contractions on a
Hilbert space $H$ and let $x\in H$. Then the following assertions hold.
\begin{itemize}
\item[(i)] $\lim_{t\to\infty}T(t)x= 0$ weakly if and only
if $\lim_{t\to\infty} \langle T(t)x,x \rangle = 0$.
\item[(ii)] If $(T(t))_{t\geq 0}$ is completely non-unitary, i.e., if there is no reducing subspace on which it is unitary, then $(T(t))_{t\geq 0}$ is weakly stable.
\end{itemize}
\end{prop}

\begin{prop}
Let $(T(t))_{t\geq0}$ be a $C_0$-semigroup of contractions on a
Hilbert space $H$. Then $(T(t)_{\mid_{W^{\perp}}})_{t\geq 0}$ has no
weakly stable orbit, hence the spectral measures of its generator
are non-Rajchman.
\end{prop}
\noindent (For the definition of Rajchman measures and a brief
discussion see Example \ref{example:Rajchman}.)

We now turn to sufficient conditions for weak stability proved
partly in Chill, Tomilov \cite{chill/tomilov:2003}. It is based on
the behaviour of the resolvent $R(\cdot, A)$ of the generator and
uses the \emph{pseudo-spectral bound of $A$} (also called
\emph{abscissa of uniform boundedness of the resolvent})
\begin{equation*}s_0(A):= \inf\bigl\{a\in \RR: R(\lambda,A) \text{ is bounded on } \{\lambda: \Re\lambda >a\}\bigr\}. \end{equation*}
\begin{thm} \label{thm:weakstab:integral}
Let $(T(t))_{t \geq 0}$ be a $C_0$-semigroup on a Banach space $X$
with generator $A$ satisfying $s_0(A) \leq 0$. Further, let $x\in X$
and $y\in X'$ be fixed. Consider the following assertions.
\begin{itemize}
\item[(a)]\hfil
$\displaystyle \int_0^1 \int_{-\infty}^{\infty} |\langle
R^2(a+is,A)x,y \rangle| \dd s\dd a < \infty$.
\item[(b)]\hfil$\displaystyle\lim_{a\to 0+}a \int_{-\infty}^{\infty} |\langle R^2(a+is,A)x,y
\rangle| \dd s = 0$.
\item[(c)] \hfil$\displaystyle \lim_{t\to \infty} \langle T(t)x,y \rangle = 0$
\end{itemize}
Then $(a)\Rightarrow (b)\Rightarrow (c)$. In particular, if (a) or
(b) holds for all $x \in X$ and $y \in X'$, then $(T(t))_{t \geq 0}$
is weakly stable.
\end{thm}
\begin{proof}
First we show that (a) implies (b).

Assume that (a) holds. From the theory of Hardy spaces we know that the function \\
$f:(0,1)\mapsto \RR_+$ defined by
\begin{equation*}
f(a):=\int_{-\infty}^{\infty} | \langle R^2(a+is,A)x,y
\rangle | \dd s\end{equation*} is monotone decreasing for $a>0$ (see Rosenblum, Rovnyak
\cite{rosenblum/rovnyak:1994} for the theory of Hardy spaces). 
%
Assume now that (b) is not true. Then there exists a monotonic decreasing null sequence $\{a_n\}_{n=1}^\infty$
such that 
\begin{equation}\label{eq:af(a)}
a_n
f(a_n)
\geq c \end{equation}
holds for some $c>0$ and all $n\in N$. 

Take now any $n,m\in \NN$ such that $a_n\leq \frac{a_m}{2}$. 
By (\ref{eq:af(a)}) and monotonicity of $f$ we have
%
\begin{equation*}
\int_{a_n}^{a_m}f(a)da \geq \sum_{k=n}^{m-1}(a_k-a_{k+1})f(a_k)
\geq \frac{c}{a_n} (a_m - a_n)=c \left(\frac{a_m}{a_n}-1\right)\geq c 
\end{equation*}
holds. This contradicts $(a)$ and the implication $(a)\Rightarrow (b)$ is proved.
%

It remains to show that $(b)$ implies $(c)$.

By (b) we have for every $a>0$
\begin{equation*}
\int_{-\infty}^{\infty} |\langle R^2(a+is,A)x,y \rangle| \dd s <
\infty. \end{equation*}
Moreover, condition $s_0(A)\leq 0$ implies that the function
$\lambda \mapsto \langle R^2(\lambda,A)x,y \rangle$ is bounded on
every half-plane $\{\lambda:\Re\, \lambda \geq a\}$. Therefore, it
belongs to the Hardy space $H^1(\{\lambda:\Re \lambda > a \})$
and
\begin{equation*}\int_{-\infty}^{\infty} |\langle R^2(a+is,A)x,y \rangle| \dd s < \infty \end{equation*}
holds for all $a>0$.
This allows us to represent the semigroup as the inverse Laplace
transform for all $a>\max\{0, \omega_0(T)\}$, where $\omega_0(T)$ is
the growth bound of $(T(t))_{t\geq 0}$. Indeed, from e.g. Kaashoek,
Verduyn Lunel \cite{kaashoek/lunel:1994} or Kaiser, Weis
\cite{kaiser/weis:2003} it follows that
\begin{equation} \label{eq:weakstab:intrep}
\langle T(t)x,y \rangle = \frac{1}{2\pi t}  \int_{-\infty}^{\infty}
\ee^{(a+is)t} \langle R^2(a+is,A)x,y \rangle \dd s.
\end{equation}
A standard application of Cauchy's theorem extends the validity of
\eqref{eq:weakstab:intrep} to all $a>0$. We now take $t =
\frac{1}{a}$ to obtain
\begin{equation*}\left| \langle T(t)x,y \rangle \right|
  \leq a \int_{-\infty}^{\infty} | \langle R^2(a+is,A)x,y \rangle | \dd s \to 0\end{equation*}
as $a \to 0+$, so $t=\frac{1}{a} \to \infty$.
\end{proof}

The implication $(b)\Rightarrow (c)$ is stated in Chill, Tomilov
\cite{chill/tomilov:2003}. They also show that the strong analogues
of $(a)$ and $(b)$ both imply strong stability of the semigroup.
Note that the relation $(a) \Rightarrow (b)$ is also valid for the
strong case by the same arguments.

\smallskip

We conclude this section with the following remarkable fact about
weak stability. By Theorem \ref{str.stable_1} one has almost weak
stability under quite general assumptions. As we will see in the
next section, almost weak stability does not imply weak
stability and, moreover, the difference between these two concepts
is fundamental (see Theorem \ref{thm:category}). In particular, this
means that weak convergence of the semigroup to zero along some
sequence $\{t_n\}_{n=1}^\infty$ with $t_n\to \infty$ in general does
not imply weak stability. However, once the sequence
$\{t_n\}_{n=1}^\infty$ is relatively dense, 
i.e., there exists a
number $\ell>0$ such that every sub-interval of $\RR_+$ of length $\ell$
intersects $\{t_n:n\in\NN\}$ (see Bart, Goldberg \cite{bart/goldberg:1978} for the terminology),  one \emph{does obtain} weak
stability.

\begin{thm} \label{thm:equidistant}
Let $(T(t))_{t \geq 0}$ be a $C_0$-semigroup on a Banach space $X$.
Suppose that \\ $\lim\limits_{n\to\infty} T(t_n)= 0$ weakly for some
relatively dense sequence $\{t_n\}_{n=1}^\infty \subset \RR_+$. Then
$(T(t))_{t \geq 0}$ is weakly stable.
\end{thm}
\begin{proof}Without loss of generality assume that $\{t_n\}_{n=1}^\infty$ is
monotone increasing and set $\ell:=\sup_{n\in \NN}(t_{n+1}-t_n)$,
which is finite by assumption. Since every $C_0$-semigroup is
bounded on compact time intervals, and $(T(t_n))_{n\in\NN}$ is
weakly converging, hence bounded, we obtain that the semigroup
$(T(t))_{t\geq 0}$ is bounded.

Fix $x\in X$, $y\in X'$. For $t\in [t_{n}, t_{n+1}]$ we have
\begin{equation*}
\langle T(t)x,y \rangle = \langle T(t-t_n)x,T'(t_n)y \rangle,
\end{equation*}
where $(T'(t))_{t\geq 0}$ is the adjoint semigroup. We note that by
assumption $T'(t_n)y\to 0$ in the weak*-topology.

Further, the set $K_x:=\{T(s)x:\ 0\leq s\leq \ell\}$ is compact in
$X$ and $T(t-t_n)x\in K_x$ for every $n\in\NN$. Since pointwise
convergence is equivalent to the uniform convergence on compact sets
(see, e.g., \cite{engel/nagel:2000}, Prop.~A.3), we see that
$\langle T(t)x,y \rangle \to 0$.
\end{proof}

Note that taking $t_n:=n$ in Theorem \ref{thm:equidistant} we obtain
that $(T(t))_{t\geq 0}$ is weakly stable if and only if $T(n)\to 0$
weakly as $n\to\infty$, $n\in\NN$. This gives a connection between
weak stability of discrete and continuous semigroups.

\section{Examples} \label{section:examples}

In this section we discuss concrete and abstract examples of almost
weakly but not weakly stable\- semigroups. Finally, we present recent
results showing that weakly and almost weakly stable semigroups even
have different Baire category in spaces of unitary and isometric
operators on Hilbert spaces.

The first example indicates how one can construct almost weakly but
not weakly stable semigroups using dynamical systems arising in
ergodic theory.

\begin{example} \label{ergodic theory}
A measurable measure-preserving semiflow $(\varphi_t)_{t \geq 0}$ on
a probability space \linebreak $(\Omega,\mathcal M,\mu)$ is called
\emph{strongly mixing} if $\displaystyle \lim_{t\to\infty} \mu(
\varphi_t^{-1}(A) \cap B) = \mu(A)\mu(B)$ for any two measurable
sets $A,B \in \mathcal M$. The semiflow $(\varphi_t)_{t \geq 0}$ is
called \emph{weakly mixing} if for all $A,B \in \mathcal M$ we have
\begin{equation*} \lim_{t\to \infty} \frac{1}{t} \int_0^t |\mu(\varphi_s^{-1}(A)
\cap B) - \mu(A)\mu(B)| \dd s = 0.\end{equation*}
These concepts play an essential role in ergodic theory, and we
refer to the monographs Cornfeld, Fomin, Sinai
\cite{cornfeld/fomin/sinai:1982}, Krengel \cite{krengel:1985},
Petersen \cite{petersen:1983}, or Halmos \cite{halmos:1956} for
further information. Clearly, strong mixing implies weak mixing, but
the converse implication does not hold in general. However, examples
of weakly but not strongly mixing semiflows are not easy to
construct; see Lind \cite{lind75} for an example and Petersen
\cite{petersen:1983}, p.~209 for a method of constructing such
semiflows.

The semiflow $(\varphi_t)_{t\geq 0}$ on $(\Omega,\mathcal M,\mu)$
induces a semigroup of isometries $(T(t))_{t\geq 0}$ on each of the
Banach spaces $X=L^p(\Omega,\mu)$ $(1\leq p < \infty)$ by defining
\begin{equation*}
(T(t)f)(\omega) := f(\varphi_t(\omega)),\ \ \omega \in \Omega, f \in
L^p(\Omega,\mu).
\end{equation*}
This semigroup is strongly continuous (see Krengel
\cite{krengel:1985}, {\S}1.6, Thm.~6.13) and relatively weakly compact
by virtue of Example \ref{ex:relcomp} (b) with $u=\Eins$. It is
well-known (see, e.g., Halmos \cite{halmos:1956}, pp.~37--38) that
\begin{equation*}
(\varphi_t)_{t\geq 0} \text{ is strongly mixing  }
 \Longleftrightarrow
\lim_{t\to\infty} \langle T(t)f,g \rangle = \langle P f,g \rangle
\text{  for all  } f\in X, \ g\in X',
\end{equation*}
and
\begin{equation*}
(\varphi_t)_{t\geq 0} \text{ is weakly mixing} \Longleftrightarrow
 \lim_{t\to\infty} \frac{1}{t}\int_0^t |\langle
T(s)f,g\rangle - \langle Pf,g\rangle|\dd s =0 \text{   for all }
f\in X,\ g\in X',
\end{equation*}
where $P$ is the projection onto $\Fix(T)$ given by $Pf:=\int_\Omega
f \dd\mu\cdot \Eins$ for all $f\in X$. Note that in both cases
$\Fix(T)=\langle \Eins \rangle$ holds.

Take now any semiflow $(\varphi_t)_{t\geq 0}$ which is weakly but
not strongly mixing. Observe that $X=X_0\oplus \langle \Eins
\rangle$, where
\begin{equation*}
X_0:= \biggl\{ f \in X \,:\, \int_{\Omega} f \dd\mu = 0 \biggr\}
\end{equation*}
is closed and $(T(t))_{t\geq 0}$-invariant.
We denote the restriction of $(T(t))_{t \geq 0}$ to $X_0$ by
$(T_0(t))_{t\geq 0}$ and its generator  by $A_0$. The semigroup
$(T_0(t))_{t\geq 0}$ is still relatively weakly compact and, since
$\Psigma (A)\cap i \RR =\emptyset$, it is almost weakly stable. On
the other hand, $(T_0(t))_{t\geq 0}$ is not weakly stable since
$(\varphi_t)_{t\geq 0}$ is not strongly mixing.

\hspace{0.1cm}

We can also look at this example from a different perspective. If
$(\varphi_t)_{t \in \RR}$ is even a measure preserving \emph{flow}, it
induces a $C_0$-group $(T(t))_{t \in \RR}$ of unitary operators on
the Hilbert space $L^2(\Omega,\mu)$. Hence we can apply the spectral
theorem and obtain for each $x\in H$ a measure $\nu_x$ on $\RR$ such
that
\begin{equation*}
\langle T(t)x,x \rangle = \int_{\RR} \ee^{i tr} \dd \nu_x(r) \ \ \
\text{for all } t\geq 0.
\end{equation*}
Thus $\langle T(t)x,x \rangle$ becomes the Fourier transform of the
measure $\nu_x$. In the next example we classify these measures
according to the behaviour of their Fourier transform at infinity.
\end{example}

\begin{example} \label{example:Rajchman}
Let us consider the Hilbert space $H=L^2(\RR,\mu)$, where $\mu$ is a
finite positive Borel measure, and the operator $A$ on $H$ is the
multiplication operator
\begin{equation*}
Af(r) = i rf(r),\ \ r \in \RR,
\end{equation*}
on its maximal domain. Then $A$ generates the unitary group
$(T(t)f)(r):=\ee^{i tr}f(r)$. Since Hilbert spaces are reflexive,
$(T(t))_{t\geq 0}$ is relatively weakly compact by Example
\ref{ex:relcomp} (a).

Clearly, $\sigma(A) \subseteq i \RR$ and $i r \in i \RR$ is an
eigenvalue of $A$ if and only if $\mu(\{r\})
> 0$. Hence, if $\mu( \{r\} ) = 0$ for all $r \in \RR$, then $A$ has
no eigenvalues and the Jacobs-Glicksberg-de Leeuw decomposition
yields that $(T(t))_{t \geq 0}$ is almost weakly stable.

For $f,g \in H$ we have  $\langle T(t)f,g \rangle = \int_\RR \ee^{i
tr}f(r)\overline{g}(r) \dd \mu$. In particular, by taking $f=g=\Eins
$ we obtain $\langle T(t)\Eins,\Eins \rangle = \int_\RR \ee^{i tr}
\dd \mu = \mathcal F \mu (t)$, the Fourier transform of $\mu$. On
the other hand, $\lim_{t\to\infty}\mathcal F \mu (t)= 0$ implies
$\lim_{t\to\infty}\langle T(t)f,g \rangle= 0$ for all $f,g \in H$,
therefore
\begin{equation*}
(T(t))_{t\geq 0} \text{  is weakly stable  } \Longleftrightarrow \
\lim_{t\to\infty} \mathcal F \mu (t)= 0.
\end{equation*}
Note that since for unitary groups weak stability as $t\to \infty$
coincides with weak stability as $t\to -\infty$, the property above
is equivalent to
\begin{equation*}\lim_{|t|\to\infty} \mathcal F \mu (t)= 0.\end{equation*}

In harmonic analysis, this property of the measure $\mu$ got its own
name. Indeed,  $\mu$ is called \emph{Rajchman} if its Fourier transform
vanishes at infinity. We refer to Lyons \cite{lyons85, lyons95} for
a brief historical overview on these measures and their properties.

We note that absolutely continuous measures are always Rajchman by
the Riemann-Lebesgue lemma and all Rajchman measures are continuous
by Wiener's theorem. However, there are continuous measures which
are not Rajchman and Rajchman measures which are not absolutely
continuous (see Lyons \cite{lyons95}).
It is
now a consequence of the considerations above that each continuous
non-Rajchman measure gives rise to an almost weakly but not weakly
stable unitary group. In Engel, Nagel \cite{engel/nagel:2000}, p.~316 
an example of a unitary group even with bounded generator  
is given, for which the corresponding spectral measures are not Rajchman.
%
\end{example}

Next, we give an example of a positive semigroup on a Banach lattice
which is almost weakly stable but not weakly stable.

\begin{example} \label{fluss}
As in Nagel (ed.) \cite{nagel:1986}, p.~206, we start from a flow on
$\CC \backslash \{0\}$ with the following properties:
\begin{enumerate}[1)]
\item The orbits starting in $z$ with $|z|\neq 1$ spiral towards the unit circle $\Gamma$;
\item $1$ is the fixed point of $\varphi$ and $\Gamma \setminus \{1\}$ is a homoclinic orbit, i.e., $\displaystyle \lim_{t\to -\infty}\varphi_t(z) = \displaystyle \lim_{t\to \infty}\varphi_t(z) =1$ for every $z\in \Gamma$.
\end{enumerate}
A concrete example comes from the differential equation in polar
coordinates $(r,\omega)=(r(t), \omega(t))$:
\begin{equation*}\left\{
\begin{array}{rcl}
\dot{r}&=&1-r, \\
\dot{\omega}&=& 1+(r^2 - 2r \cos\omega).
\end{array}\right.
\end{equation*}
Take $x_0\in \CC$ with $0<|x_0|<1$ and denote by
$S_{x_0}:=\{\varphi_t(x_0):\:t\geq 0\}$ the orbit starting from
$x_0$. Then $S:=S_{x_0}\cup \Gamma$ is compact for the usual
topology of $\CC$.

\smallskip We define a multiplication on $S$ as follows. For
$x=\varphi_t(x_0)$ and $y=\varphi_s(x_0)$ we put
\begin{equation*}
x y := \varphi_{t+s}(x_0).
\end{equation*}
For $x\in \Gamma$, $x=\lim_{n\to\infty} x_n$,
$x_n=\varphi_{t_n}(x_0)\in S_{x_0}$ and $y=\varphi_s(x_0)\in
S_{x_0}$, we define $x y =y x := \lim_{n\to\infty} x_n y$.
Note that by $|x_n y - \varphi_s(x)|=|\varphi_s
(x_n)-\varphi_s(x)|\leq C |x_n - x| \underset{n\to
\infty}{\longrightarrow} 0$ the definition is correct and satisfies
\begin{equation*}
x y = \varphi_s(x).
\end{equation*}
For $x,y\in \Gamma$ we define $xy:=1$. This multiplication on $S$ is
separately continuous and makes $S$ a semi-topological semigroup
(see Engel, Nagel \cite{engel/nagel:2000}, Sec.~V.2).

Consider now the Banach space $X:=C(S)$. By Example \ref{ex:relcomp}
(\ref{ex:semitop}) the set
\begin{equation*}
\{ f(s\:\cdot):\: s\in S \} \subset C(S)
\end{equation*}
is relatively weakly compact for every $f\in C(S)$. By definition of
the multiplication on $S$ this implies that
\begin{equation*}
\{ f(\varphi_t(\cdot)):\: t\geq 0 \}
\end{equation*}
is relatively weakly compact in $C(S)$. Consider the semigroup
induced by the flow, i.e.,
\begin{equation*}
(T(t)f)(x):=f(\varphi_t(x)),\quad\ f\in C(S),\: x\in S.
\end{equation*}
By the above, each orbit $\{T(t)f:\: t\geq 0 \}$ is relatively
weakly compact in $C(S)$ and hence, by Lemma \ref{weakcomp},
$(T(t))_{t\geq 0}$ is weakly compact. Note that the strong
continuity of $(T(t))_{t\geq 0}$ follows, as shown in Nagel (ed.)
\cite{nagel:1986}, Lemma B-II.3.2, from the separate continuity of
the flow. Furthermore, the semigroup $(T(t))_{t\geq 0}$ is
isometric.

Next, we take $X_0:=\{f\in C(S):\: f(1)=0\}$ and identify it with
the Banach lattice $C_0(S\setminus\{1\})$. Then both subspaces in
the decomposition $C(S)=X_0\oplus \langle \Eins  \rangle$ are
invariant under $(T(t))_{t\geq 0}$. Denote by $(T_0(t))_{t\geq 0}$
the restricted semigroup to $X_0$ and by $A_0$ its generator. The
semigroup $(T_0(t))_{t\geq 0}$ is still relatively weakly compact.

Since $\Fix(T_0):=\bigcap_{t\geq 0}\Fix(T_0(t))=\{0\}$, we have that
$0\notin \Psigma(A_0)$. Moreover, $\Psigma (A_0) \cap i\RR =
\emptyset$ holds, which implies by the Jacobs-Glicksberg-de Leeuw
theorem that $(T_0(t))_{t\geq 0}$ is almost weakly stable.

To see that $(T_0(t))_{t\geq 0}$ is not weakly stable it is enough
to consider $\delta_{x_0}\in X_0'$. Since
\begin{equation*}
\langle T_0(t)f, \delta_{x_0} \rangle = f(\varphi(t,x_0)), \ f\in
X_0,
\end{equation*}
$f(\Gamma)$ always belongs to the closure of $\{\langle T_0(t)f,
\delta_{x_0} \rangle:\: t\geq 0\}$ and hence the semigroup
$(T_0(t))_{t\geq 0}$ can not be weakly stable.
\end{example}

\noindent Let us summarise the above as follows.
\begin{prop}\label{prop:cexa1}There exist a locally compact space $\Omega$ and a positive, relatively weakly compact $C_0$-semigroup
of isometries on $C_0(\Omega)$ which is almost weakly but not weakly
stable.
\end{prop}

This 
enables us to answer a question of Emelyanov \cite{emelyanov:2005}
in negative. Consider the discrete semigroup $(T(n))_{n\in\NN}=(T(1)^n)_{n\in\NN}$ 
from Proposition \ref{prop:cexa1}. 
By a result of Jones and Lin \cite{jones/lin:1980}, we know that $0$ belongs 
to the weak closure of each of the orbits. Whereas Theorem \ref{thm:equidistant} 
shows that this semigroup is not weakly stable. The semigroup is positive and 
isometric on the Banach lattice $C_0(\Omega)$. 

 Moreover, Proposition \ref{prop:cexa1} becomes particularly
interesting in view of the following results; for details and
discussion see Chill, Tomilov \cite{chill/tomilov:2004}.

\begin{thm}[{Groh, Neubrander \cite{groh/neubrander:1981}, Thm.~3.2}; {Chill, Tomilov \cite{chill/tomilov:2004}, Thm.~7.7}]
\label{thm:AL,AM} For a bounded, positive, mean ergodic
$C_0$-semigroup $(T(t))_{t\geq 0}$ on a Banach lattice $X$
 with generator $(A,D(A))$, the following assertions hold.
\begin{enumerate}[(i)]
\item
If $X\cong L^1(\Omega, \mu)$, then $\Psigma(A)\cap i\RR = \emptyset$
is equivalent to the \emph{strong stability} of $(T(t))_{t\geq 0}$.
\item \label{AM}
If $X\cong C(K)$, K compact, then $\Psigma(A)\cap i\RR = \emptyset$
is equivalent to the \emph{uniform exponential stability} of
$(T(t))_{t\geq 0}$.
\end{enumerate}
\end{thm}
Example \ref{fluss} shows that we can not drop the assumption on the
existence of a unit element in $X$ in Theorem \ref{thm:AL,AM}
(\ref{AM}).

\hspace{0.5cm}

At the end of this section, we show that the examples presented
above represent the general situation. Indeed, typical isometric semigroups and typical unitary groups are almost
weakly but not weakly stable in the following sense. 
%
\begin{thm}\label{thm:category}
Let $H$ be a separable infinite-dimensional Hilbert space and let
$\mathcal{I}$ denote the set of all isometric $C_0$-semigroups on
$H$. Then the set of all weakly stable isometric semigroups is of
first category and the set of all almost weakly stable isometric
semigroups is residual in $\mathcal{I}$ with respect to an 
appropriate topology making  $\mathcal{I}$ to a completely metrisable space. 
The analogous statement holds for unitary groups.
\end{thm}
\noindent For the proof and precise description of the appropriate
topologies on the spaces of all unitary and all isometric semigroups 
see Eisner, Ser\'eny \cite{eisner/sereny-cont:2006}, and
\cite{eisner/sereny:2006} for the case of unitary, isometric and
contractive operators. These results are the operator theoretic
counterpart of classical category theorems of Halmos and Rokhlin
from ergodic theory, see Halmos \cite{halmos:1956}, pp.~77--80.

\section{Individual stability and local resolvent} \label{section:individual_stability}

In this section, we restrict our attention to single orbits and
present results implying
\begin{equation*}
\lim_{t\to\infty} \langle T(t)x, y \rangle =0 \quad \text{ for some given  } x \text {  and } y. 
\end{equation*}
%
%
The tool will be the \emph{bounded local resolvent} $R(\lambda)x_0$
which exists by definition if the function $\rho(A)\ni\lambda\mapsto
R(\lambda,A)x_0$ admits a bounded, holomorphic extension
$R(\lambda)x_0$ to the whole right half-plane $\{\lambda:\Re \lambda
>0\}$. This we assume in the following.

Clearly, if we suppose that for \emph{all} $x_0\in X$ the local
resolvent $R(\lambda)x_0$ is bounded,  analyticity and the principle
of uniform boundedness yield the boundedness of the operator
resolvent $R(\lambda, A)$ on $\{\lambda:\Re \lambda
>0\}$, hence uniform exponential
stability for semigroups on Hilbert spaces and (at least) strong
stability for semigroups on reflexive Banach spaces (see the
theorems of Gearhart and Arendt, Batty, Lyubich, V\~u, e.g., in
Engel, Nagel \cite{engel/nagel:2000}, Thm.~V.1.11 and V.2.21). The
reasonable questions therefore address the individual stability of a
single orbit in terms of the local resolvent of one single element
$x_0\in X $.

Without any differentiability or boundedness assumption on the
semigroup it is necessary to do some initial smoothing on $x_0$ in
order to have stability in any sense.
 However, if the Banach space $X$
has some nice geometric properties, even
strong stability can be derived. Huang and van Neerven
\cite{huang/vanneerven:1999} proved that if the Banach space is
$B$-convex or has the analytic Radon-Nikod\'ym property, then the
existence of a bounded local resolvent $R(\lambda)x_0$ on
$\{\lambda:\Re \lambda
>0\}$ already implies  \emph{strong convergence} $T(t)R(\mu,A)^\alpha
x_0\lm 0$ as $t\lm+\infty$ for any $\alpha>1$. (Here $\mu$ is
greater  than the growth bound $\omega_0(A)$, thus $R(\mu,A)$ is a
sectorial operator admitting fractional powers.) Actually, if $X$
has Fourier type $p>1$, then we can take $\alpha>1/p$, see
\cite{huang/vanneerven:1999} (the assumption that $\alpha>1/p$ is
shown to be optimal by Wrobel \cite{wrobel:1999}), and if the
semigroup is eventually differentiable and $p=2$, then no smoothing
is needed, i.e., $\alpha\geq 0$ is allowed, see \cite{huang:1999}.

In general, without any additional assumptions on the space or on
the regularity of the semigroup one can only deduce weak individual
stability. The following result is due to Huang and van Neerven
\cite{huang/vanneerven:1999}.
\begin{thm} \label{thm:weakstab:locres}
Let $(T(t))_{t \geq 0}$ be a $C_0$-semigroup on a Banach space $X$
with generator $A$. Let $x_0 \in X$ and suppose that the function
$\lambda \mapsto R(\lambda,A)x_0$ has a bounded holomorphic
extension to $\{\lambda:\Re \lambda >0\}$. Then
\begin{itemize}
\item[\it a)] $\lim_{t\to \infty} \langle T(t)x_0,y \rangle = 0$ for
all $y \in D((A')^2)$;
\item[\it b)]  $\lim_{t\to\infty} T(t)x_0=0$  weakly, if additionally the semigroup $(T(t))_{t
\geq 0}$ is uniformly bounded.
\end{itemize}
\end{thm}
Tauberian theorems are among the primary tools to deduce information
on the asymptotic behaviour of the semigroup from properties of the
resolvent, and they have extensively been used to obtain strong and
weak stability.
 We refer the reader
to the monograph Arendt, Batty, Hieber, Neubrander
\cite{arendt/etal:2001} and also to Chill \cite{chill:1998}. As
illustration we include here a proof for part b) which uses Ingham's
Tauberian theorem.
\begin{thm}[{Ingham\rm}]
Let $f : \RR_+ \to \CC$ be bounded and uniformly continuous and
suppose that the Laplace transform $\hat{f}$ of $f$ has a locally
integrable boundary function on the imaginary axis (that is, there
exists $h \in L^1_{\mathrm{loc}}(\RR,\CC)$ such that $\lim_{a\to 0+}
\hat{f}(a+i\cdot) = h$ in the distributional sense).
 Then $\lim_{t\to \infty}
f(t) =0$.
\end{thm}
\noindent For proofs and a detailed treatment see
\cite{arendt/etal:2001}, Sect.~4 and \cite{chill:1998}.
\begin{proof}[Proof of part b).]
 By assumption the operator resolvent $R(\lambda, A)x_0$ and the local resolvent $R(\lambda) x_0$ coincide on the right halfplane. So for a fix $y \in
 X'$, on the right half-plane  $\lambda\mapsto \langle R(\lambda)x_0,y\rangle$ is the Laplace transform of
 the function $t\mapsto \langle T(t)x_0,y\rangle$.
 Since $(T(t))_{t \geq 0}$
is uniformly bounded, the weak orbit $t \mapsto \langle T(t)x_0,y
\rangle$ is bounded and uniformly continuous, so we can apply
Ingham's theorem to obtain $\lim_{t\to\infty} \langle T(t)x_0,y
\rangle = 0$.
\end{proof}
For the proof of Theorem \ref{thm:weakstab:locres} part a) one could
also use Ingham's Theorem, and check the assumptions of the theorem
by following the lines of Batty, Chill and van Neerven in
\cite{batty/chill/vanneerven:2000}. Actually in
\cite{batty/chill/vanneerven:2000}  a powerful functional calculus
 method is developed, which among other yields the proof of the more general Theorem \ref{thm:weakstab:beta} below.
 To prove part a) we nevertheless choose
a different, fairly elementary way
 (see Eisner, Farkas \cite{eisner/farkas:2006}). 

\begin{proof}[Proof of part a).] By $\lambda \mapsto R(\lambda)x_0$ we denote the holomorphic
continuation of $\lambda \mapsto R(\lambda,A)x_0$ to the half-plane
$\{\lambda:\Re \lambda >0 \}$. The uniqueness theorem for
holomorphic functions and the resolvent identity imply
\begin{multline*}
R(\delta+is)x_0=R(a+i s,A)x_0+(a-\delta)
R(a+i s,A)R(\delta+i s)x_0 \\
=R(a+i s,A)x_0+(a-\delta)R^2(a+i s,A)x_0 +(a-\delta)^2 R^2(a+i
s,A)R(\delta+i s)x_0.
\end{multline*}
For all $y\in D({A'}^2)$ we have
\begin{multline*}
2\pi\ee^{-\delta t} \la T(t)x_0,y\ra=\int_{-\infty}^{\infty} \ee^{i
st}\la R(\delta+i s)x_0,y\ra\dd s =\int_{-\infty}^{\infty} \ee^{i
st}\la R(a+i s,A)x_0,y\ra\dd
s\\
+(a-\delta)\int_{-\infty}^{\infty} \ee^{i st} \la R^2(a+i
s,A)x_0,y\ra\dd s +(a-\delta)^2\int_{-\infty}^{\infty} \ee^{i st}
\la R^2(a+i s,A)R(\delta+i s)x_0,y\ra\dd s.
\end{multline*}
Indeed, for $a>\omega_0(T)$ the first equality follows from
representing the semigroup as the inverse Laplace transform of the
resolvent; subsequently, the Cauchy theorem implies this
representation for all $a>0$. The functions $f_\delta(s):=\la
R^2(a+i s,A)R(\delta+i s)x_0,y\ra$ form a relatively compact subset
of $L^1(\RR)$, because
\begin{multline*}
|f_\delta(s)|=|\la R^2(a+i s,A)R(\delta+i s)x_0,y\ra|=\\ = |\la
R(\delta+i s)x_0,R^2(a+i s,A')y\ra| \leq M \|R^2(a+i s,A')y\|,
\end{multline*}
and the function on the right hand side lies in $L^1(\RR)$, so the
family $f_\delta$ is uniformly integrable (and bounded), thus
relatively compact. By compactness we find a sequence $\delta_n\lm
0$ such that $\lim_{n\to\infty} f_{\delta_n} = f$ in $L^1(\RR)$.
 By substituting $\delta_n$ in the above equality
and letting $n\lm\infty$ we obtain
\begin{multline*}
2\pi \la T(t)x_0,y\ra=\int_{-\infty}^{\infty} \ee^{i st}\la R(a+i
s,A)x_0,y\ra\dd s
\\ +a\int_{-\infty}^{\infty} \ee^{i st} \la R^2(a+i
s,A)x_0,y\ra\dd s + a^2\int_{-\infty}^{\infty} \ee^{i st} f(s)\dd s
=  I_1(t)+I_2(t)+I_3(t).
\end{multline*}
It is easy to deal with the last term $I_3$. The function $f$ lies
in $L^1(\RR)$, so by the Riemann-Lebesgue Lemma its Fourier
transform vanishes at infinity, i.e., $\lim_{t\to\infty} I_3(t)= 0$.
Since $y\in D((A')^2)$,  we can integrate by parts in $I_1$ to
obtain
\begin{align*}
I_1(t)=\int_{-\infty}^{\infty} \ee^{i st}\la x_0,R(a+i s,A')y\ra\dd
s=\frac{1}{t}\int_{-\infty}^{\infty} \ee^{i st}\la x_0,R^2(a+i
s,A')y\ra\dd s.
\end{align*}
The last integral is absolutely convergent, because using the
resolvent identity for $R(\lambda, A')$, one can show that for $y\in
D({A'}^2)$ and $a>\omega_0(T)$ fixed, $\|R^2(a+i
s,A')y\|=O(({a^2+s^2})^{-1})$ holds. Hence
\begin{equation*}
|I_1(t)|\leq\frac{1}{t}\int_{-\infty}^{\infty} \|
x_0\|\cdot\|R^2(a+i s,A')y\|\dd s\lm 0\qquad\qquad\mbox{as
$t\lm\infty$}.
\end{equation*}
Concerning $I_2$ we observe that  $\la x_0,R^2(a+i \cdot,A')y\ra\in
L^1(\RR)$, so once again by the Riemann-Lebesgue Lemma we have
\begin{equation*}
I_2(t)=a\int_{-\infty}^{\infty} \ee^{i st} \la x_0,R^2(a+i
s,A')y\ra\dd s \lm 0\qquad\qquad\mbox{as $t\lm\infty$},
\end{equation*}
and the proof is complete.
\end{proof}

\smallskip
\noindent Actually, Huang and van Neerven
 proved the following more general
theorem.

\begin{thm}[{Huang, van Neerven \cite{huang/vanneerven:1999}}] \label{thm:weakstab:beta}
 Let $(T(t))_{t \geq 0}$ be a $C_0$-semigroup on a Banach
space $X$ with generator $A$. For $x_0 \in X$ assume that the
bounded local resolvent exists. Then $\lim_{t\to\infty}
T(t)(\lambda_0-A)^{-\beta}x_0 =0$ weakly for all $\beta>1$ and
$\lambda_0 > \omega_0(T)$.
\end{thm}

Under a special positivity condition one can take $\beta=1$ in
Theorem \ref{thm:weakstab:beta}.
\begin{thm}[{van Neerven \cite{vanneerven:2002}}]\label{thm:one}
Suppose that $X$ is an ordered Banach space with weakly closed
normal cone $C$. If for some $x_0\in X$ the function $\lambda
\mapsto R(\lambda,A)x_0$ has a bounded holomorphic extension to
$\{\lambda:\Re\lambda > 0\}$ and $T(t)x_0\in C$ for all sufficiently
large $t$, then $\lim_{t\to\infty} T(t)R(\mu,A)x_0 =0$ weakly for
all $\mu\in\rho(A)$.
\end{thm}
\noindent The above eventual positivity assumption cannot be
omitted. Indeed, van Neerven \cite{vanneerven:2002} proved that the
existence of a bounded local resolvent $R(\lambda)x_0$ in general implies
$\|T(t)R(\mu, A)x_0\|=O(1+t)$, and Batty
\cite{batty:2003} showed this to be optimal, whereas weak
convergence of $T(t)R(\mu, A)x_0$ to zero would imply $\|T(t)R(\mu,
A)x_0\|=O(1)$.

\section{Comments and open questions} \label{section:comments_and_open_quetions}
 In this closing section, we collect some
noteworthy open questions. We also touch upon two further areas
 connected to the topic of this survey: weakly almost periodic
functions and the cogenerator of contraction semigroups on Hilbert
spaces.

\smallskip
We start by listing some problems arising from this paper.
\begin{question}
In Remark \ref{rem:strongfulint} we have seen that in Hilbert spaces
strong stability is fully characterised by an integrability
condition of the resolvent of the generator. This raises the
question if the same is true for weak stability, i.e., for example, does the
converse
 of Theorem \ref{thm:weakstab:integral} hold in Hilbert spaces? Note the
difference between Theorem \ref{thm:weakstab:integral} b) and
Theorem \ref{str.stable_1} (v), the latter being equivalent to
almost weak stability.
\end{question}

\begin{question}
Understanding weak stability of unitary groups is of special
importance, but a satisfactory description is still lacking. We
mention here that not even the boundedness of the generator of the
group would make the problem easier (see Engel and Nagel
\cite{engel/nagel:2000}, p.~316).

\end{question}

\begin{question}
 By Theorem \ref{thm:equidistant} it suffices to understand weak stability of
discrete semigroups. However, a description of those sequences
$\{t_n\}_{n=1}^\infty\subseteq\RR$ for which
$\lim_{n\to\infty}T(t_n)= 0$ weakly implies $\lim_{t\to\infty} T(t)=
0$ weakly is desirable.
\end{question}

\subsection{Weakly almost periodic functions}\label{subsec:weakper}

In this part, we sketch very briefly some results of Ruess and
Summers (\cite{ruess/summers:1986}--\cite{ruess/summers:1992jmaa}) in connection with weak asymptotic behaviour of operator
semigroups. We select only those aspects of their theory that are
directly related to weak stability of orbits. However, to illustrate
the merit of their approach, we shall consider the more general
setting of \emph{almost orbits}, too. In the sequel, we assume that
$(T(t))_{t\geq 0}$ is a uniformly bounded $C_0$-semigroup with
generator $A$. A function $u:\RR_+\to X$ is called an \emph{almost
orbit} of the semigroup $(T(t))_{t\geq 0}$ if
\begin{equation*}
\lim_{t\to\infty}\sup_{h\in\RR_+}\|u(t+h)-T(h)u(t)\|=0\quad\mbox{holds}.
\end{equation*}
The use of this notion is explained for instance by the fact that a
solution of an inhomogeneous abstract Cauchy problem is an almost
orbit of the corresponding semigroup.

We shall need the following definition; see
\cite{ruess/summers:1986} and the references therein for details and
historical remarks (cf also Bart, Goldberg \cite{bart/goldberg:1978}).

\begin{defi}A bounded, continuous function $f\in C_b(\RR_+,X)$ is called
\emph{asymptotically almost periodic} if its translates
$H(f):=\{f(\cdot\:+h):h\in \RR_+\}$ form a relatively compact set in
$C_b(\RR_+,X)$ for the sup-norm. A function $f\in C_b(\RR_+,X)$ is
called \emph{Eberlein-weakly almost periodic}, if $H(f)$ is
relatively weakly compact in $C_b(\RR_+,X)$. The spaces of such
functions are denoted by $AAP(\RR_+,X)$ and $W(\RR_+,X)$
respectively, and $W_0(\RR_+,X)$ is the vector subspace of functions
$f\in W(\RR_+,X)$ such that the weak closure of $H(f)$ contains $0$.
A function $f\in C_b(\RR_+,X)$ is \emph{weakly asymptotically almost
periodic} if for all $y\in X'$ the function $y\circ f:\RR_+\to\CC$
belongs to $AAP(\RR_+,\CC)$.
\end{defi}
Clearly, an asymptotically almost periodic function is
Eberlein-weakly almost periodic and weakly asymptotically almost
periodic, while the converse implications fail in general.

It is not surprising that the function $t\mapsto T(t)x$ is
Eberlein-weakly almost periodic if and only if the orbit
$\{T(t)x:t\geq 0\}$ is relatively weakly compact (this is also true
for almost orbits, see \cite{ruess/summers:1990pac}, Sect.~4.5). The
structure theory for Eberlein-weakly almost periodic functions
developed by Ruess and Summers thus gives a slightly finer
description of relatively weakly compact semigroups than the
Glicksberg-Jacobs-de Leeuw decomposition. Using this approach Ruess
and Summers proved the following (for orbits this is basically a
corollary of the Jacobs-Glicksberg-de Leeuw decomposition, see
Theorem \ref{str.stable_1}).
\begin{thm}[{Ruess, Summers \cite{ruess/summers:1990pac}}] Let
$(T(t))_{t\geq 0}$ be a uniformly bounded $C_0$-semigroup with
generator $A$ on the reflexive Banach space $X$. Then the following are
equivalent:
\begin{enumerate}[a)]
\item $\Psigma(A)\cap i\RR=\emptyset$.
\item Every almost orbit of $(T(t))_{t\geq 0}$ belongs to
$W_0(\RR_+,X)$.
\item Every  orbit of $(T(t))_{t\geq 0}$ belongs to
$W_0(\RR_+,X)$.
\item $\liminf_{t\to\infty} |\langle
T(t)x,y\rangle|=0$ for each $x\in X$, $y\in X'$.
\end{enumerate}
\end{thm}

As for weak stability, asymptotic almost periodicity is of great use.
\begin{thm}[{Ruess, Summers \cite{ruess/summers:1992jmaa}}]
 Let $u$ be a weakly asymptotically almost periodic almost orbit of $(T(t))_{t\geq 0}$ and suppose that its range $\{u(t):t\geq 0\}$ is relatively weakly compact.
  Then the following
 hold:
\begin{enumerate}[a)]
\item If $\Psigma(A)\cap i\RR=\emptyset$, then $\lim_{t\to\infty}
u(t)= 0$ weakly.
\item If $\Psigma(A)\cap i\RR=\{0\}$, then $u(t)$ has a weak limit,
  which is
a fixed point of $(T(t))_{t\geq 0}$.
\end{enumerate}
\end{thm}

For a detailed treatment, as well as for connections to ergodic
theory and further generalisations we refer to
\cite{ruess/summers:1986,ruess/summers:1988,ruess/summers:1990a,ruess/summers:1990pac,ruess/summers:1992pam, ruess/summers:1992jmaa}
(see also Arendt, Batty, Hieber, Neubrander \cite{arendt/etal:2001},
Section~5.4).

\subsection{Cogenerator of contractive semigroups}
Let $(T(t))_{t\geq 0}$ be a contractive $C_0$-semigroup on a Hilbert
space $H$. The \emph{cogenerator} of $(T(t))_{t\geq 0}$ is defined
as the (negative) Cayley-transform of the infinitesimal generator
$A$ of $(T(t))_{t\geq 0}$, i.e.,
\begin{equation*}
G:=-(I+A)(I-A)^{-1}=I-2R(1,A).
\end{equation*}
It is easy to see that the cogenerator is a contraction, see
Sz.-Nagy, Foia{\c{s}} \cite{sznagy/foias}, Sections~III.8--9 for
details. Using a functional calculus the semigroup can be directly
recovered and, moreover, many important properties of the semigroup
can be read off from its cogenerator. Namely, the semigroup consists
of normal, self-adjoint, isometric or unitary operators if and only
if the cogenerator is normal, self-adjoint, isometric or unitary,
respectively. Also the asymptotic behaviour of the semigroup can be
characterised with the help of the cogenerator.
\begin{thm}\label{thm:cogen}
Let $(T(t))_{t\geq 0}$ be a $C_0$-semigroup of contractions and $G$
its cogenerator. Then
\begin{equation*}
\lim_{t\to\infty}\|T(t)x\|=\lim_{n\to\infty}\|G^nx\|.
\end{equation*}
In particular, the semigroup is strongly stable if and only if $G$
is strongly stable.
\end{thm}
Motivated by this we ask the following.
\begin{question} Is the analogue of Theorem \ref{thm:cogen} true for weak stability?
More precisely, is there a connection between the weak
stability of a contractive semigroup and the weak stability of its
cogenerator?
\end{question}
Note that the function $z\mapsto-\frac{1+z}{1-z}$ maps the imaginary
axis onto the unit circle, so by the spectral mapping theorem for
the point spectrum (see Engel, Nagel \cite{engel/nagel:2000}, Thm.~IV.3.7), we
have that
\begin{equation*}
\Psigma(A)\cap i\RR=\emptyset\quad\mbox{if and only if}\quad \Psigma
(G)\cap \{z:|z|=1\}=\emptyset.
\end{equation*}
Hence by a result of Jones and Lin \cite{jones/lin:1980} we obtain
the following.
\begin{prop}
A contraction semigroup $(T(t))_{t\geq0}$ on a Hilbert space $H$ is
almost weakly stable if and only if its cogenerator $G$ is ``almost
weakly stable'', i.e., when $0$ belongs to the weak closure of each
orbit $\{G^nx:n\in\NN\}$, $x\in H$.
\end{prop}

This again connects the asymptotic behaviour of $(T(t))_{t\geq0}$ 
to the behaviour of the powers of a single operator, thus allowing 
the application of results by, e.g., Jones and Lin \cite{jones/lin:1976, jones/lin:1980}.


\parindent0pt
\end{document}